\documentclass[12pt,reqno]{amsart}
\newtheorem{thm}{Theorem}[section]
\newtheorem{defi}{Definition}[section]

\newtheorem{cor}{Corollary}[section]
\newtheorem{pr}{Proposition}[section]
\theoremstyle{definition}
\newtheorem{rem}{Remark}[section]
\newtheorem{exm}{Example}[section]
\newcommand{\be}{\begin{equation}}
\newcommand{\ee}{\end{equation}}
\newcommand{\bea}{\begin{eqnarray}}
\newcommand{\eea}{\end{eqnarray}}
\newcommand{\beb}{\begin{eqnarray*}}
\newcommand{\eeb}{\end{eqnarray*}}
\usepackage{amssymb,amsfonts,amsthm,setspace,indentfirst,mathrsfs}
\usepackage{minibox}
\usepackage{enumitem}
\numberwithin{equation}{section}
\begin{document}
%%%%%%%%%%%%%%%%%%%%%%%%%%%%%%%%%%%%%%%%%%%%%%%%%%%%%%%%%%%%%%%%%%%%%%%%%%%%%%%%%%%%%%%%%%%%%%%%%%%%%%%%%%%%%
%
\title[Curvature properties of generalized pp-wave metric]{\bf{Curvature properties of generalized pp-wave metric}}
\author[A. A. Shaikh, T. Q. Binh and H. Kundu]{Absos Ali Shaikh$^{1}$, Tran Quoc Binh$^2$ and Haradhan Kundu$^3$}
\date{\today}
\address{\noindent$^{1}$ Department of Mathematics,
\newline Aligarh Muslim University,
\newline Aligarh-202002,
\newline Uttar Pradesh, India}
\email{aask2003@yahoo.co.in, aashaikh@math.buruniv.ac.in}
\address{\noindent$^{2}$ Institute of Mathematics and Informatics,
\newline Debrecen University,
\newline H-4010 Debrecen, P.O. Box 12, Hungary}
\email{binh@math.klte.hu}
\address{\noindent$^3$ Department of Mathematics,
\newline University of Burdwan, Golapbag,
\newline Burdwan-713104,
\newline West Bengal, India}
\email{kundu.haradhan@gmail.com}
\dedicatory{}
%%%%%%%%%%%%%%%%%%%%%%%%%%%%%%%%%%%%%%%%%%%%% Abstract %%%%%%%%%%%%%%%%%%%%%%%%%%%%%%%%%%%%%%%%%%%%%%%%%%%%%
\begin{abstract}
The main objective of the present paper is to investigate the curvature properties of generalized pp-wave metric. It is shown that generalized pp-wave spacetime is Ricci generalized pseudosymmetric, 2-quasi-Einstein and generalized quasi-Einstein in the sense of Chaki. As a special case it is shown that pp-wave spacetime is semisymmetric, semisymmetric due to conformal and projective curvature tensors, $R$-space by Venzi and satisfies the pseudosymmetric type condition $P\cdot P = -\frac{1}{3}Q(S, P)$. Again we investigate the sufficient condition for which a generalized pp-wave spacetime turns into pp-wave spacetime, pure radiation spacetime, locally symmetric and recurrent. Finally, it is shown that the energy-momentum tensor of pp-wave spacetime is parallel if and only if it is cyclic parallel. And the energy momentum tensor is Codazzi type if it is cyclic parallel but the converse is not true as shown by an example. Finally we make a comparison between the curvature properties of the Robinson-Trautman metric and generalized pp-wave metric.
\end{abstract}
%%%%%%%%%%%%%%%%%%%%%%%%%%%%%%%%%%%%%%%%%%%%%%%%%%%%%%%%%%%%%%%%%%%%%%%%%%%%%%%%%%%%%%%%%%%%%%%%%%%%%%%%%%%%
%
\maketitle
%
%%%%%%%%%%%%%%%%%%%%%%%%%%%%%%%%%%%%%%%%%%%%%%%%%%%%%%%%%%%%%%%%%%%%%%%%%%%%%%%%%%%%%%%%%%%%%%%%%%%%%%%%%%%%%%
%																					Introduction
%%%%%%%%%%%%%%%%%%%%%%%%%%%%%%%%%%%%%%%%%%%%%%%%%%%%%%%%%%%%%%%%%%%%%%%%%%%%%%%%%%%%%%%%%%%%%%%%%%%%%%%%%%%%%%
\section{\bf Introduction}\label{intro}
%%%%%%%%%%%%%%%%%%%%%%%%%%%%%%%%%%%%%%%
A spacetime is a 4-dimensional connected Lorentzian manifold. The class of pp-wave metrics (\cite{EK62}, \cite{ste03}) arose during the study of exact solutions of Einstein's field equations. The term ``pp-wave'' is given by Ehlers and Kundt \cite{EK62}, where ``pp'' stands for ``plane-fronted gravitational waves with parallel rays'', which means that the generalized pp-wave spacetimes admit geodesic null vector field whose twist, expansion and shear are zero. The ``plane rays'' implies that the rotation of the vector field vanishes. For vacuum type N, this implies the existence of a covariantly constant vector field which is parallel to the null vector field. There are various forms of generalized pp-wave metrics in different coordinates. The pp-wave belongs to the class of solutions admitting a non-expanding, shear-free and twist-free null congruence and it admits a null Killing vector.\\
%===================================================================
\indent The family of pp-wave space-times was first discussed by Brinkmann \cite{Brin25} and interpreted in terms of gravitational waves by Peres \cite{Pere59}. According to Brinkmann, a pp-wave spacetime is any Lorentzian manifold
whose metric tensor can be described, with respect to Brinkmann coordinates, in the form
\be\label{pp-bc}
ds^2 = H(u, x, y) du^2 + 2 du dv + dx^2 + dy^2,
\ee
where $H$ is any nowhere vanishing smooth function. Again it is well known that a Lorentzian manifold with parallel lightlike (null) vector field is called Brinkmann-wave (\cite{Brin25}, \cite{GL10}). A Brinkmann-wave is called pp-wave if its curvature tensor $R$ satisfies the condition ${R_{ij}}^{pq}R_{pqkl}= 0$ (\cite{GL10}, \cite{MS16}, \cite{SG86}). In 1984, Radhakrishna and Singh \cite{RS84} presented a class of solutions to Einstein-Maxwell equation for the null electrovac Petrov type $N$ gravitational field. They presented a metric of the form
\beb
ds^2=-2U du^2+2dudr-\frac{1}{2}P[(dx^3)^2+(dx^4)^2],
\eeb
where $U=U(u,x^3,x^4)$ and $P=P(x^3,x^4)$ are two nowhere vanishing smooth functions. 
For the simplicity of notation, we write the variable $u$ as $x$, and the function $U$ as $h$ and $P$ as $f$. Then the aforesaid metric can be written as
\be\label{gppwm}
ds^2=-2 h(x,x^3,x^4) (dx)^2 + 2 dx dr-\frac{1}{2}f(x^3,x^4)[(dx^3)^2+(dx^4)^2].
\ee
%========================================================================
\indent In section \ref{se-gpp} we show that the metric \eqref{gppwm} admits a covariantly constant null vector field and it satisfies the condition ${R_{ij}}^{pq}R_{pqkl}= 0$ if
\be\label{cond}
\left(f f_{33} - f^2_{3}\right)+\left(f f_{44} - f^2_{4}\right) = 0.
\ee
Thus we can say that the metric \eqref{gppwm} is a Brinkmann-wave and it becomes a pp-wave if \eqref{cond} holds. Hence we can say the metric \eqref{gppwm} as ``generalized pp-wave metric''. We note that for $f\equiv -2$ and $h = -\frac{1}{2}H(u,x^3,x^4)$, the solution \eqref{gppwm} reduces to the pp-wave metric \eqref{pp-bc}.\\
%=======================================================================================
\indent In the study of differential geometry the notion of manifold of constant curvature has been generalized by many authors in different directions such as locally symmetric manifold by Cartan \cite{Cart26}, semisymmetric manifold by Cartan \cite{Cart46} (see also \cite{Szab82}, \cite{Szab84}, \cite{Szab85}), pseudosymmetric manifold by Adam\'{o}w and Deszcz \cite{AD83}, recurrent manifold by Ruse (\cite{Ruse46}, \cite{Ruse49a}, \cite{Ruse49b}, see also \cite{Walk50}), weakly generalized recurrent manifold by Shaikh and Roy (\cite{SAR13}, \cite{SR11}), hyper generalized recurrent manifold by Shaikh and Patra (\cite{SP10}, \cite{SRK15}), super generalized recurrent manifold by Shaikh et al. \cite{SRK16}. We note that in such geometric structures the first order and second order covariant derivatives of the Riemann curvature tensor and other curvature tensors involved. We mention that the notion of pseudosymmetry by Deszcz is important in the study of differential geometry due to its application in relativity and cosmology (see \cite{DDVV94}, \cite{DHV04}, \cite{DVV91} and also references therein). The notion of pseudosymmetry is extended by considering other curvature tensors in its defining condition, such as conformal pseudosymmetry, pseudosymmetric Weyl conformal curvature tensor, Ricci generalized pseudosymmetry etc. and they are called pseudosymmetric type conditions. It may be mentioned that different pseudosymmetric type conditions are admitted by various spacetimes, such as G\"{o}del spacetime (\cite{DHJKS14}, \cite{Gode49}), Som-Raychaudhuri spacetime (\cite{SK16srs}, \cite{SR68}), Reissner-Nordstr\"{o}m spacetime \cite{Kowa06} and Robertson-Walker spacetime (\cite{ADEHM14}, \cite{DDHKS00}, \cite{DK99}).\\
%===========================================================================
\indent The main object of the present paper is to investigate the geometric structures admitted by the generalized pp-wave metric \eqref{gppwm}. It is interesting to note that the metric \eqref{gppwm} without any other condition admits several geometric structures, such as  Ricci generalized pseudosymmetry, 2-quasi Einstein and generalized quasi-Einstein in the sense of Chaki \cite{Chak01}. Again it is shown that the pp-wave metric (i.e., \eqref{gppwm} with condition \eqref{cond}) is Ricci recurrent but not recurrent, semisymmetric, $R$-space by Venzi, conformal curvature 2-forms are recurrent, Ricci tensor is Riemann compatible and fulfills a pseudosymmetric type condition due to the projective curvature tensor $P$. For the study of pseudosymmetric type conditions with projective curvature tensor we refer the reader to see the recent papers of Shaikh and Kundu (\cite{SKppsn}, \cite{SKppsnw}). It is interesting to note that for such a metric $P\cdot R = 0$ but $P\cdot\mathcal R \ne 0$.\\
%=============================================================================
\indent It is also shown that the metric is weakly Ricci symmetric and weakly cyclic Ricci symmetric for different associated 1-forms, which ensures the existence of infinitely many solutions of associated 1-forms of such structures. Again we investigate the condition for which such a spacetime is locally symmetric and recurrent.\\
%=============================================================================
\indent The paper is organized as follows. Section 2 deals with defining conditions of different curvature restricted geometric structures, such as recurrent, semisymmetry, pseudosymmetry, weakly symmetry etc. as preliminaries. Section 3 is devoted to the investigation of curvature restricted geometric structures admitted by the generalized pp-wave metric \eqref{gppwm}. Section 4 is mainly concerned with the geometric structures admitted by pp-wave metric and plane wave metric. Section 5 deals with the investigation of the conditions under which the energy-momentum tensor of such spacetimes are parallel, Codazzi type and cyclic parallel. Finally, the last section is devoted to make a comparison between the curvature properties of the Robinson-Trautman metric and generalized pp-wave metric as well as pp-wave metric.
%%%%%%%%%%%%%%%%%%%%%%%%%%%%%%%%%%%%%%%%%%%%%%%%%%%%%%%%%%%%%%%%%%%%%%%%%%%%%%%%%%%%%%%%%%%%%%%%%%%%%%%%%%%%%%%%%%%%
%                                                     Preliminaries
%%%%%%%%%%%%%%%%%%%%%%%%%%%%%%%%%%%%%%%%%%%%%%%%%%%%%%%%%%%%%%%%%%%%%%%%%%%%%%%%%%%%%%%%%%%%%%%%%%%%%%%%%%%%%%%%%%%%
\section{\bf Preliminaries}
%%%%%%%%%%%%%%%%%%%%%%%%%%%
Let $M$ be a connected smooth semi-Riemannian manifold of dimension $n$ $(\geq 3)$ equipped with the semi-Riemannian metric $g$. Let $R$, $\mathcal R$, $S$, $\mathcal S$ and $\kappa $ be respectively the Riemann-Christoffel curvature tensor of type $(0,4)$, the Riemann-Christoffel curvature tensor of type $(1,3)$, the Ricci tensor of type $(0,2)$, the Ricci tensor of type $(1,1)$ and the scalar curvature of $M$.\\
%=================================================================
\indent For two symmetric $(0,2)$-tensors $A$ and $E$, their Kulkarni-Nomizu product $A\wedge E$ is defined as (see e.g. \cite{DGHS11}, \cite{Glog02}):
\begin{eqnarray*}
(A\wedge E)(X_1,X_2,X_3,X_4) &=& A(X_1,X_4)E(X_2,X_3) + A(X_2,X_3)E(X_1,X_4)\\
&-& A(X_1,X_3)E(X_2,X_4) - A(X_2,X_4)E(X_1,X_3),
\end{eqnarray*}
where $X_1,X_2,X_3,X_4 \in \chi(M)$, the Lie algebra of all smooth vector fields on $M$. Throughout the paper we will consider $X, Y, X_1, X_2, \cdots \in \chi(M)$.\\
%===========================================================
\indent In terms of Kulkarni-Nomizu product, the conformal curvature tensor $C$, the concircular curvature tensor $W$, the conharmonic curvature tensor $K$ (\cite{Ishi57}, \cite{YK89}) and the Gaussian curvature tensor $\mathfrak G$ can be expressed as
\beb
C &=& R-\frac{1}{n-2}(g\wedge S) + \frac{r}{2(n-2)(n-1)}(g\wedge g),\\
W &=& R-\frac{r}{2n(n-1)}(g\wedge g),\\
K &=& R-\frac{1}{n-2}(g\wedge S),\\
\mathfrak G &=& \frac{1}{2}(g\wedge g).
\eeb
Again the projective curvature tensor $P$ of type $(0,4)$ is given by
$$
P(X_1, X_2, X_3, X_4) = R(X_1, X_2, X_3, X_4) - \frac{1}{n-1}[g(X_1, X_4)S(X_2, X_3)-g(X_2, X_4)S(X_1, X_3)].
$$
%=============================================
\indent For a symmetric $(0, 2)$-tensor $A$, we get an endomorphism $\mathcal A$ defined by $g(\mathcal AX,Y) = A(X,Y)$. Then its $k$-th level tensor $A^k$ of type $(0,2)$ is given by
$$A^k(X,Y) = A(\mathcal A^{k-1}X,Y),$$
where $\mathcal A^{k-1}$ is the endomorphism corresponding to $A^{k-1}$.\\
%===========================================================
\begin{defi} (\cite{Bess87}, \cite{SKgrt})
A semi-Riemannian manifold $M$ is said to be $Ein(2)$, $Ein(3)$ and $Ein(4)$ respectively if
$$S^2 + \lambda_1 S + \lambda_2 g = 0,$$
$$S^3 + \lambda_3 S^2 + \lambda_4 S + \lambda_5 g = 0 \ \mbox{and}$$
$$S^4 + \lambda_6 S^3 + \lambda_7 S^2 + \lambda_8 S + \lambda_9 g = 0$$
holds for some scalars $\lambda_i, \, 1\le i \le 9$.
\end{defi}
%=================================================
\begin{defi} (\cite{Bess87}, \cite{SKgrt})
A semi-Riemannian manifold $M$ is said to be generalized Roter type (\cite{SKgrt}, \cite{SK16}) if its Riemann-Christoffel curvature tensor $R$ can be expressed as a linear combination of $g\wedge g$, $g\wedge S$, $S\wedge S$, $g\wedge S^2$, $S\wedge S^2$ and $S^2\wedge S^2$. Again $M$ is said to be Roter type (\cite{Desz03}, \cite{Desz03a}) if $R$ can be expressed as a linear combination of $g\wedge g$, $g\wedge S$ and $S\wedge S$.
\end{defi}
%=====================================
\indent Again for a $(0,4)$-tensor $D$, an endomorphism $\mathscr{D}(X,Y)$ and the corresponding $(1,3)$-tensor $\mathcal D$ can be defined as
$$\mathscr{D}(X,Y)X_1 = \mathcal D(X,Y)X_1, \ \ D(X_1,X_2,X_3,X_4) = g(\mathcal D(X_1,X_2)X_3, X_4).$$
Again for a symmetric $(0,2)$-tensor $A$, another endomorphism $X\wedge_A Y$ (\cite{DDHKS00}, \cite{DGHS11}) can be defined as
$$(X\wedge_A Y)X_1 = A(Y,X_1)X-A(X,X_1)Y.$$
%================================================
\indent By operating $\mathscr{D}(X,Y)$ and $X\wedge_A Y$ on a $(0,k)$-tensor $B$, $k\geq 1$, we can obtain two $(0,k+2)$-tensors $D\cdot B$ and $Q(A,B)$ respectively given by (see \cite{DG02}, \cite{DGHS98}, \cite{DH03}, \cite{SDHJK15}, \cite{SK14} and also references therein):
$$D\cdot B(X_1,X_2,\cdots,X_k,X,Y) = -B(\mathcal D(X,Y)X_1,X_2,\cdots,X_k) - \cdots - B(X_1,X_2,\cdots,\mathcal D(X,Y)X_k).$$
\beb
\mbox{and } \ &&Q(A,B)(X_1,X_2, \ldots ,X_k,X,Y) = ((X \wedge_A Y)\cdot B)(X_1,X_2, \ldots ,X_k)\\
&&= A(X, X_1) B(Y,X_2,\cdots,X_k) + \cdots + A(X, X_k) B(X_1,X_2,\cdots,Y)\\
&& - A(Y, X_1) B(X,X_2,\cdots,X_k) - \cdots - A(Y, X_k) B(X_1,X_2,\cdots,X).
\eeb
In terms of local coordinates system, $D\cdot B$ and $Q(A,B)$ can be written as
\beb
(D\cdot B)_{i_1 i_2\cdots i_k j l} &=& -g^{pq}\left[D_{jli_1 q}B_{p i_2\cdots i_k} + \cdots + D_{jli_k q}B_{i_1 i_2\cdots p}\right],
\eeb
\beb
Q(A, B)_{i_1 i_2\cdots i_k j l} &=& A_{li_1}B_{j i_2\cdots i_k} + \cdots + A_{li_k}B_{i_1 i_2\cdots j}\\
																	& & -A_{ji_1}B_{l i_2\cdots i_k} - \cdots - A_{ji_k}B_{i_1 i_2\cdots l}.
\eeb
%===============================================
\begin{defi}$($\cite{Ruse46}, \cite{Ruse49a}, \cite{Ruse49b}, \cite{Patt52}$)$
A semi-Riemannian manifold $M$ is said to be $B$-recurrent if $\nabla B = \Pi \otimes B$ for an 1-form $\Pi$. In particular for $B = R$ (resp., $S$), the manifold $M$ is called recurrent (resp., Ricci recurrent) manifold.
\end{defi}
%=========================================
\begin{defi}$($\cite{AD83}, \cite{Cart46}, \cite{Desz92}, \cite{SK14}, \cite{Szab82}$)$
A semi-Riemannian manifold $M$ is said to be $B$-semisymmetric type if $D\cdot B = 0$ and it is said to be $B$-pseudosymmetric type if $\left(\sum\limits_{i=1}^k c_i D_i\right)\cdot B = 0$ for some scalars $c_i$'s, where $D$ and each $D_i$, $i=1,\ldots, k$, $(k\ge 2)$, are (0,4) curvature tensors. In particular, if $c_i$'s are all constants, then it is called $B$-pseudosymmetric type manifold of constant type or otherwise non-constant type.
\end{defi}
%=======
\indent In particular, if $i =2$, $D_1 = R$, $D_2 = \mathfrak G$ and $B = R$, then $M$ is called Ricci generalized pseudosymmetric \cite{DD91a}.\\
%===============================================
\begin{defi} 
A semi-Riemannian manifold $M$ is said to be quasi-Einstein (resp., 2-quasi-Einstein) if at each point of $M$, rank$(S - \alpha g)\le 1$ (resp., $\le 2$) for a scalar $\alpha$. Also $M$ is said to be generalized quasi-Einstein in the sense of Chaki \cite{Chak01} if
$$S = \alpha g + \beta \Pi \otimes \Pi + \gamma [\Pi \otimes \Omega + \Pi \otimes \Omega]$$
form some 1-forms $\Pi$ and $\Omega$.
\end{defi}
%
%================================================================================
\begin{defi}\label{def2.8}
Let $D$ be a $(0,4)$-tensor and $E$ be a symmetric $(0, 2)$-tensor on $M$. Then $E$ is said to be $D$-compatible (\cite{DGJPZ13}, \cite{MM12b}, \cite{MM13}) if
\[
D(\mathcal E X_1, X,X_2,X_3) + D(\mathcal E X_2, X,X_3,X_1) + D(\mathcal E X_3, X,X_1,X_2) = 0
\]
holds, where $\mathcal E$ is the endomorphism corresponding to $E$ defined as $g(\mathcal E X_1, X_2) = E(X_1, X_2)$. Again an 1-form $\Pi$ is said to be $D$-compatible if $\Pi\otimes \Pi$ is $D$-compatible.
\end{defi}
%=================================================================================
\begin{defi}
A Riemannian manifold $M$ is said to be weakly cyclic Ricci symmetric \cite{SJ06} if its Ricci tensor satisfies the condition
\bea\label{crs}
&&(\nabla_X S)(X_1,X_2) + (\nabla_{X_1} S)(X,X_2) + (\nabla_{X_2} S)(X_1,X)\\\nonumber
&&=\Pi(X)\, S(X_1,X_2) + \Omega(X_1)\, S(X,X_2) + \Theta(X_2)\, S(X_1,X),
\eea
for three 1-forms $\Pi$, $\Omega$ and $\Theta$ on $M$. Such a manifold is called weakly cyclic Ricci symmetric manifold with solution $(\Pi, \Omega, \Theta)$. Moreover if the first term of left hand side is equal to the right hand side, then it is called weakly Ricci symmetric manifold \cite{TB93}.
\end{defi}
%======================================================================
\begin{defi}
Let $D$ be a $(0,4)$ tensor and $Z$ be a $(0,2)$-tensor on $M$. Then the corresponding curvature 2-forms $\Omega_{(D)l}^m$ (\cite{Bess87}, \cite{LR89}) are called recurrent if and only if (\cite{MS12a}, \cite{MS13a}, \cite{MS14})
\beb\label{man}
&&(\nabla_{X_1} D)(X_2,X_3,X,Y)+(\nabla_{X_2} D)(X_3,X_1,X,Y)+(\nabla_{X_3} D)(X_1,X_2,X,Y) =\\
&&\hspace{1in} \Pi(X_1) D(X_2,X_3,X,Y) + \Pi(X_2) D(X_3,X_1,X,Y)+ \Pi(X_3) D(X_1,X_2,X,Y)
\eeb
and 1-forms $\Lambda_{(Z)l}$ \cite{SKP03} are called recurrent if and only if
$$(\nabla_{X_1} Z)(X_2,X) - (\nabla_{X_2} Z)(X_1,X) = \Pi(X_1) Z(X_2,X) - \Pi(X_2) Z(X_1,X)$$
for an 1-form $\Pi$.
\end{defi}
%-------------------------------------------------------------------
\begin{defi}$($\cite{Prav95}, \cite{SKppsn}, \cite{Venz85}$)$ 
Let $\mathcal L(M)$ be the vector space formed by all 1-forms $\Theta$ on $M$ satisfying
$$\Theta(X_1)D(X_2,X_3,X_4,X_5)+\Theta(X_2)D(X_3,X_1,X_4,X_5)+\Theta(X_3)D(X_1,X_2,X_4,X_5) = 0,$$
where $D$ is a $(0,4)$-tensor. Then $M$ is said to be a $D$-space by Venzi if $dim \mathcal L(M) \ge 1$.
\end{defi}
%%%%%%%%%%%%%%%%%%%%%%%%%%%%%%%%%%%%%%%%%%%%%%%%%%%%%%%%%%%%%%%%%%%%%%%%%%%%%%%%%%%%%%%%%%%%%%%%%%%%%%%%%%%%%%%%%%%%%%%%%
%                                                  Generalized pp-wave metric
%%%%%%%%%%%%%%%%%%%%%%%%%%%%%%%%%%%%%%%%%%%%%%%%%%%%%%%%%%%%%%%%%%%%%%%%%%%%%%%%%%%%%%%%%%%%%%%%%%%%%%%%%%%%%%%%%%%%%%%%%
\section{\bf Curvature properties of generalized pp-wave metric}\label{se-gpp}
%%%%%%%%%%%%%%%%%%%%%%%%%%%%%%%%%%%%%%%%%%%%%%%%%%%%%%%%%%%%%%%%%%%%%%%%%%%%%%
We can now write the metric tensor $g$ of the generalized pp-wave metric \eqref{gppwm} as follows:
$$g = \left(
\begin{array}{cccc}
 -2 h & 1 & 0 & 0 \\
 1 & 0 & 0 & 0 \\
 0 & 0 & -\frac{1}{2} f & 0 \\
 0 & 0 & 0 & -\frac{1}{2} f
\end{array}
\right).$$
%=========================================================================
Then the non-zero components of its Riemann-Christoffel curvature tensor $R$, Ricci tensor $S$ and scalar curvature $\kappa$ of \eqref{gppwm} are given by
$$R_{1313}=\frac{-f_3 h_3+f_4 h_4+2 f h_{33}}{2 f}, \ \ R_{1314}=\frac{-f_4 h_3-f_3 h_4+2 f h_{34}}{2 f},$$
$$R_{1414}=\frac{f_3 h_3-f_4 h_4+2 f h_{44}}{2 f}, \ \ R_{3434}=\frac{-f_3^2-f_4^2+f f_{33}+f f_{44}}{4 f},$$
$$S_{11}=\frac{2 \left(h_{33}+h_{44}\right)}{f}, \ \ S_{33}= S_{44}=\frac{-f_3^2-f_4^2+f f_{33}+f f_{44}}{2 f^2},$$
$$\kappa = \frac{2 \left(f_3^2+f_4^2-f \left(f_{33}+f_{44}\right)\right)}{f^3}.$$
%-------------------------------------------------------------------------------------------------------------
Again the non-zero components of $C$ and $P$ are given by
$$C_{1313}= -C_{1414}=\frac{-f_3 h_3+f_4 h_4+f h_{33}-f h_{44}}{2 f}, \ \ C_{1314}=\frac{-f_4 h_3-f_3 h_4+2 f h_{34}}{2 f},$$
$$P_{1211}=\frac{2 \left(h_{33}+h_{44}\right)}{3 f}, \ P_{1313}=\frac{-3 f_3 h_3+3 f_4 h_4+4 f h_{33}-2 f h_{44}}{6 f},$$
$$P_{1314}= -P_{1341}= P_{1413}= -P_{1431}=\frac{-f_4 h_3-f_3 h_4+2 f h_{34}}{2 f}, \ \ P_{1441}=-\frac{f_3 h_3-f_4 h_4+2 f h_{44}}{2 f},$$
$$P_{1331}=-\frac{-f_3 h_3+f_4 h_4+2 f h_{33}}{2 f}, \ \ P_{1414}=-\frac{-3 f_3 h_3+3 f_4 h_4+2 f h_{33}-4 f h_{44}}{6 f}.$$
%-----------------------------------------------------------------------------------------------------
Now the non-zero components (upto symmetry) of the energy momentum tensor 
$$T= \frac{c^4}{8\pi G}\left[S-\left(\frac{\kappa}{2}-\Lambda\right)g\right],$$ 
where $c=$ speed of light in vacuum, $G=$ gravitational constant and $\Lambda=$ cosmological constant, are given by
$$T_{11}=-\frac{c^4 \left(f^3 h \Lambda -f^2 h_{33}-f^2 h_{44}+f_{33} f h+f_{44} f h-f_3^2 h-f_4^2 h\right)}{4 \pi  f^3 G},$$
$$T_{12}=\frac{c^4 \left(f^3 \Lambda +f_{33} f+f_{44} f-f_3^2-f_4^2\right)}{8 \pi  f^3 G}, \ \ T_{33}= T_{44}=-\frac{c^4 f \Lambda }{16 \pi  G}.$$
Then the non-zero components (upto symmetry) of covariant derivative $\nabla T$ of the energy momentum tensor $T$ are given by
$$T_{11,1}=\frac{c^4 \left(h_{133}+h_{144}\right)}{4 \pi  f G},$$
\beb
T_{11,3}=\frac{c^4}{4 \pi  f^4 G} &&\left(-f^2 h f_{344}-f^2 h f_{333}+f^3 h_{333}+f^3 h_{344}-f_3 f^2 h_{33}-f_3 f^2 h_{44}\right.\\
&&\left.+4 f_3 f_{33} f h+2 f_4 f_{34} f h+2 f_3 f_{44} f h-3 f_3^3 h-3 f_3 f_4^2 h\right),
\eeb
\beb
T_{11,4}=\frac{c^4}{4 \pi  f^4 G} &&\left(-f^2 h f_{444}-f^2 h f_{344}+f^3 h_{334}+f^3 h_{444}-f_4 f^2 h_{33}-f_4 f^2 h_{44}\right.\\
&&\left.+2 f_4 f_{33} f h+2 f_3 f_{34} f h+4 f_4 f_{44} f h-3 f_4^3 h-3 f_3^2 f_4 h\right),
\eeb
$$T_{12,3}=\frac{c^4 \left(f^2 f_{344}+f^2 f_{333}+3 f_3^3+3 f_4^2 f_3-4 f f_{33} f_3-2 f f_{44} f_3-2 f f_4 f_{34}\right)}{8 \pi  f^4 G},$$
$$T_{12,4}=\frac{c^4 \left(f^2 f_{444}+f^2 f_{334}+3 f_4^3+3 f_3^2 f_4-2 f f_{33} f_4-4 f f_{44} f_4-2 f f_3 f_{34}\right)}{8 \pi  f^4 G},$$
$$T_{13,1}=\frac{c^4 \left(-f_3^2-f_4^2+f f_{33}+f f_{44}\right) h_3}{8 \pi  f^3 G},$$
$$T_{14,1}=\frac{c^4 \left(-f_3^2-f_4^2+f f_{33}+f f_{44}\right) h_4}{8 \pi  f^3 G}.$$
%%%%%%%%%%%%%%%%%%%%%%%%%%%%%%%%%%%%%%%%%%%%%%%%%%%%%%%%%%%%%%%%%%%%%%%%%%%%%%%%%%%%%%%%%%%%%%%%%%%%%%%%%%%%%%%%%%%%%%%%%%%%%%
\indent From above we see that the Ricci tensor $S$ of \eqref{gppwm} is of the form
\be\label{eqchqe}
S(X, Y) = \alpha g(X, Y) + \beta \eta(X)\eta(Y) + \gamma [\eta(X)\delta(Y)+\eta(Y)\delta(X)],
\ee
where $\alpha = \frac{f_3^2+f_4^2-f f_{33}-f f_{44}}{f^3}, \ \beta = 1, \ \gamma = 1, \ \eta = \{1,0,0,0\}$\\
\indent \indent and $\delta = \left\{\frac{2 f^2 \left(h_{33}+h_{44}\right)-2 \left(f_{33}+f_{44}\right) f h+2 \left(f_3^2+f_4^2\right) h - f^3}{2 f^3}, \frac{f \left(f_{33}+f_{44}-f_3^2-f_4^2\right)}{f^3},0,0\right\}.$\\
Therefore the metric \eqref{gppwm} is generalized quasi-Einstein in the sense of Chaki. Moreover $||\eta|| = 0$, $||\delta||^2 = \frac{\left(f_3^2+f_4^2-f \left(f_{33}+f_{44}\right)\right) \left(f-2 \left(h_{33}+h_{44}\right)\right)}{f^4}$, $g(\eta, \delta) = -\frac{f_3^2+f_4^2-f \left(f_{33}+f_{44}\right)}{f^3}$ and $\nabla \eta = 0$. So there exists a null covariantly constant vector field $\zeta$, where $\zeta$ is the corresponding vector field of $\eta$ (i.e., $g(\zeta, X) = \eta(X)$, $\forall \ X$). Hence we can conclude that the spacetime with the metric \eqref{gppwm} is a generalized pp-wave metric.\\
%===================================================================================================
\indent Now from the value of the components of various tensors related to \eqref{gppwm}, we can state the following:
\begin{thm}\label{maingppwm}
The generalized pp-wave metric \eqref{gppwm} possesses the following curvature restricted geometric structures:
\begin{enumerate}[label=(\roman*)]
\item Ricci generalized pseudosymmetric such that $R\cdot R = Q(S, R)$,
\item 2-quasi-Einstein, since Rank$(S - \frac{f_3^2+f_4^2-f f_{33}-f f_{44}}{f^3}g) = 2$,
\item generalized quasi-Einstein in the sense of Chaki such that \eqref{eqchqe} holds,
\item Ricci tensor is Riemann compatible as well as Weyl compatible,
\item $Ein(3)$ manifold such that $S^2 = -\frac{f_3^2+f_4^2-f \left(f_{33}+f_{44}\right)}{f^3} S^3$.
\end{enumerate}
\end{thm}
%====================================================================================================
\indent Now from the components of $R$, we see that the only non-zero (upto symmetry) the tensor $D_{ijkl} = R_{ij}^{pq}R_{pqkl}$ is given by $D_{3434} = \frac{\left(f_3^2+f_4^2-f \left(f_{33}+f_{44}\right)\right)^2}{2 f^4}$. Hence from definition, we can state the following.
\begin{thm}
The generalized pp-wave metric \eqref{gppwm} becomes a pp-wave metric if
$$f_3^2+f_4^2-f \left(f_{33}+f_{44}\right) = 0.$$
\end{thm}
%
%%%%%%%%%%%%%%%%%%%%%%%%%%%%%%%%%%%%%%%%%%%%%%%%%%%%%%%%%%%%%%%%%%%%%%%%%%%%%%%%%%%%%%%%%%%%%%%%%%%%%%%%%%%%%%%%
%%%%%%%%%%%%%%%%%%%%%%%%%%%%%%%%%%%%%%%%%%%%%%%%%%%%%%%%%%%%%%%%%%%%%%%%%%%%%%%%%%%%%%%%%%%%%%%%%%%%%%%%%%%%%%%%
%                                                  pp-wave metric
%%%%%%%%%%%%%%%%%%%%%%%%%%%%%%%%%%%%%%%%%%%%%%%%%%%%%%%%%%%%%%%%%%%%%%%%%%%%%%%%%%%%%%%%%%%%%%%%%%%%%%%%%%%%%%%%
\section{\bf Curvature properties of pp-wave and plane wave metric}\label{se-pp}
%%%%%%%%%%%%%%%%%%%%%%%%%%%%%%%%%%%%%%%%%%%%%%%%%%%%%%%%%%%%%%%%%%%%%%%%%%%%%%%%
In this section we investigate the curvature restricted geometric structures admitted by the pp-wave metric. Since under the condition \eqref{cond} the generalized pp-wave metric \eqref{gppwm} becomes a pp-wave metric, putting this condition we get the non-zero components of $R$, $S$, $C$ and $P$ of the pp-wave metric given as follows:
$$R_{1313}=\frac{-f_3 h_3+f_4 h_4+2 f h_{33}}{2 f}, \ R_{1314}=\frac{-f_4 h_3-f_3 h_4+2 f h_{34}}{2 f}, \ R_{1414}=\frac{f_3 h_3-f_4 h_4+2 f h_{44}}{2 f},$$
$$S_{11}=\frac{2 \left(h_{33}+h_{44}\right)}{f},$$
$$C_{1313}= -C_{1414}=\frac{-f_3 h_3+f_4 h_4+f h_{33}-f h_{44}}{2 f}, \ \ C_{1314}=\frac{-f_4 h_3-f_3 h_4+2 f h_{34}}{2 f},$$
$$P_{1211}=\frac{2 \left(h_{33}+h_{44}\right)}{3 f}, \ P_{1313}=\frac{-3 f_3 h_3+3 f_4 h_4+4 f h_{33}-2 f h_{44}}{6 f},$$
$$P_{1314}= -P_{1341}= P_{1413}= -P_{1431}=\frac{-f_4 h_3-f_3 h_4+2 f h_{34}}{2 f}, \ \ P_{1441}=-\frac{f_3 h_3-f_4 h_4+2 f h_{44}}{2 f},$$
$$P_{1331}=-\frac{-f_3 h_3+f_4 h_4+2 f h_{33}}{2 f}, \ \ P_{1414}=-\frac{-3 f_3 h_3+3 f_4 h_4+2 f h_{33}-4 f h_{44}}{6 f}.$$
%=====================================================================================================
Using the values of the components of $g$, $R$, $S$ and $C$ we get\\
(i) $\kappa = 0$, (ii) $R\cdot R = 0$, (iii) $R\cdot S = 0$, (iv) $Q(S, R) = 0$, (v) $R\cdot C = 0$, (vi) $C\cdot R = 0$, (vii) $C\cdot C = 0$ and (viii) $Q(S, C) =0$.\\
%=============================================================================================================
%=============================================================================================================
\indent The energy momentum tensor $T$ is given by
\be\label{tcom}
T_{11}=-\frac{c^4 \left(f h \Lambda -h_{33}-h_{44}\right)}{4 \pi  f G}, \ \ T_{12}=\frac{c^4 \Lambda }{8 \pi  G}, \ \ T_{33}=T_{44}=-\frac{c^4 f \Lambda }{16 \pi  G}.
\ee
%--------------------------
Then the non-zero components of covariant derivative of $T$ are given by
\be\label{dtcom}
\left.\begin{array}{l}
T_{11,1}=\frac{c^4 \left(h_{144}+h_{133}\right)}{4 \pi  f G},\\
T_{11,3}=\frac{c^4 \left(f h_{344}+f h_{333}-f_3 h_{33}-f_3 h_{44}\right)}{4 \pi  f^2 G},\\
T_{11,4}=\frac{c^4 \left(f h_{444}+f h_{334}-f_4 h_{33}-f_4 h_{44}\right)}{4 \pi  f^2 G}.
\end{array}\right\}
\ee
%%%%%%%%%%%%%%%%%%%%%%%%%%%%%%%%%%%%%%%%%%%%%%%%%%%%%%%%%%%%%%%%%%%%%%%%%%%%%%%%%%%%%%%%%%%%%%%%%%%%%%%%%%%
\indent From the above calculations, we can state the following:
%====================================================================================================
\begin{thm}\label{mt}
The pp-wave metric (\eqref{gppwm} with the additional condition \eqref{cond}) possesses the following curvature properties:
\begin{enumerate}[label=(\roman*)]
%---------------------------------------------------------------------------------------------------
		\item $\kappa = 0$ and hence $R=W$ and $C=K$.
		%---------------------------------------------------------------------------------------------------
    \item $R$-space and $C$-space by Venzi for $\{1,0,0,0\}$. Hence from second Bianchi identity, the curvature 2-forms $\Omega_{(R)l}^m$ are recurrent for the 1-form $\{1,0,0,0\}$.
		%---------------------------------------------------------------------------------------------------
		\item It is semisymmetric and hence Ricci semisymmetric, conformally semisymmetric and projectively semisymmetric. 
		%--------------------------------------------------------------------------------------------------
		\item If $\alpha_{_3}+\alpha_{_4}\ne 0$, then the conformal 2-forms $\Omega_{(C)l}^m$ are recurrent with 1-form of recurrency
		$$\Pi = \left\{1,0,\frac{\alpha_{_1}+\alpha_{_2}}{\alpha_{_3}+\alpha_{_4}},\frac{\alpha_{_5}+\alpha_{_6}+\alpha_{_7}}{\alpha_{_3}+\alpha_{_4}}\right\},$$
where $\alpha_{_1} = \left(f_3 h_3-f_4 h_4-f \left(h_{33}-h_{44}\right)\right) \left(f_3 \left(h_{33}+h_{44}\right)-f \left(h_{333}+h_{344}\right)\right)$,\\
$\alpha_{_2}= \left(-f_4 h_3-f_3 h_4+2 f h_{34}\right) \left(f \left(h_{334}+h_{444}\right)-f_4 \left(h_{33}+h_{44}\right)\right),$\\
$\alpha_{_3}= f^2 \left(4 h_{34}^2+\left(h_{33}-h_{44}\right){}^2\right)+f_3^2 \left(h_3^2+h_4^2\right)+f_4^2 \left(h_3^2+h_4^2\right),$\\
$\alpha_{_4}= 2 f f_4 \left(h_4 \left(h_{33}-h_{44}\right)-2 h_3 h_{34}\right)-2 f f_3 \left(2 h_4 h_{34}+h_3 \left(h_{33}-h_{44}\right)\right),$\\
$\alpha_{_5}= f^2 \left(2 h_{34} \left(h_{333}+h_{344}\right)+\left(h_{44}-h_{33}\right) \left(h_{334}+h_{444}\right)\right)+f_3^2 h_4 \left(h_{33}+h_{44}\right)+f_4^2 h_4 \left(h_{33}+h_{44}\right),$\\
$\alpha_{_6}= f f_3 \left(-2 h_{34} \left(h_{33}+h_{44}\right)-h_4 \left(h_{333}+h_{344}\right)+h_3 \left(h_{334}+h_{444}\right)\right)$,\\
$\alpha_{_7}= -f f_4 \left(-h_{33}^2+h_{44}^2+h_3 \left(h_{333}+h_{344}\right)+h_4 \left(h_{334}+h_{444}\right)\right)$.
		%---------------------------------------------------------------------------------------------------
		\item It is not recurrent but if $h_{33}+h_{44}\ne 0$, then it is  Ricci recurrent  with 1-form of recurrency
		$$\Pi = \left\{\frac{h_{144}+h_{133}}{h_{33}+h_{44}}, 0,\frac{f h_{344}+f h_{333}-f_3 h_{33}-f_3 h_{44}}{f\left(h_{33}+h_{44}\right)},\frac{f h_{444}+f h_{334}-f_4 h_{33}-f_4 h_{44}}{f\left(h_{33}+h_{44}\right)}\right\}.$$
		%----------------------------------------------------------------------------------------------------
		\item If $h_{33}+h_{44}\ne 0$, then it is weakly cyclic Ricci recurrent with solution $(\Pi, \Omega, \Theta)$, given by
		$$\Pi = \left\{\Pi_1, 0, \frac{f \left(h_{333}+h_{344}\right)-f_3 \left(h_{33}+h_{44}\right)}{f \left(h_{33}+h_{44}\right)}, \frac{f \left(h_{334}+h_{444}\right)-f_4 \left(h_{33}+h_{44}\right)}{f \left(h_{33}+h_{44}\right)}\right\},$$
		$$\Omega = \left\{\Omega_1, 0, \frac{f \left(h_{333}+h_{344}\right)-f_3 \left(h_{33}+h_{44}\right)}{f \left(h_{33}+h_{44}\right)}, \frac{f \left(h_{334}+h_{444}\right)-f_4 \left(h_{33}+h_{44}\right)}{f \left(h_{33}+h_{44}\right)}\right\},$$
		$$\Theta = \left\{\frac{3 \left(h_{133}+h_{144}\right)}{h_{33}+h_{44}}-\Pi_1-\Omega_1, 0, \frac{f \left(h_{333}+h_{344}\right)-f_3 \left(h_{33}+h_{44}\right)}{f \left(h_{33}+h_{44}\right)}, \frac{f \left(h_{334}+h_{444}\right)-f_4 \left(h_{33}+h_{44}\right)}{f \left(h_{33}+h_{44}\right)}\right\},$$
		where $\Pi_1$ and $\Omega_1$ are arbitrary scalars.
		%----------------------------------------------------------------------------------------------------
		\item Ricci simple (i.e., $S = \alpha \eta\otimes\eta$) for
		$$\alpha = \frac{2(h_{33}+h_{44})}{f} \ \ \mbox{ and } \ \ \eta = \left\{1, 0, 0, 0\right\}$$
		and hence $S\wedge S =0$ and $S^2 = 0$. Again $||\eta|| = 0$ and $\nabla\eta = 0$.
		%-----------------------------------------------------------------------------------------------------
		\item $Q(S, R) = Q(S, C) =0$ but $R$ or $C$ is not a scalar multiple of $S\wedge S$ as $S$ is of rank 1.
		%------------------------------------------------------------------------------------------------------
		\item $C\cdot R =0$ and hence $C\cdot S =0$, $C\cdot C =0$ and $C\cdot P =0$.
		%------------------------------------------------------------------------------------------------------
		\item $P\cdot R = 0$ but $P\cdot \mathcal R \ne 0$. Also but $P\cdot S =P\cdot \mathcal S =0$.
		%------------------------------------------------------------------------------------------------------
		\item Ricci tensor is Riemann compatible as well as Weyl compatible.
		%------------------------------------------------------------------------------------------------------
		\item $P\cdot P = -\frac{1}{3}Q(S,P)$.
		%------------------------------------------------------------------------------------------------------
\end{enumerate}
\end{thm}
%
%===============================================================================================================
\begin{rem}
From the value of the local components (presented in Section \ref{se-pp}) of various tensors of the pp-wave metric, we can easily conclude that the metric is
\begin{enumerate}[label=(\roman*)]
\item not conformally symmetric and hence not locally symmetric or projectively symmetric,
%------------------------------------------------
\item not conformally recurrent and hence not recurrent or not projectively recurrent,
%------------------------------------------------
\item not super generalized recurrent \cite{SRK16} and hence not hyper generalized recurrent \cite{SP10}, weakly generalized recurrent \cite{SR11},
%------------------------------------------------
\item not weakly symmetric \cite{TB89} for $R, C, P, W$ and $K$ and hence not Chaki pseudosymmetric \cite{Chak87} for $R, C, P, W$ and $K$,
%------------------------------------------------
\item neither cyclic Ricci parallel \cite{Gray78} nor of  Codazzi type Ricci tensor although its scalar curvature is constant,
%-------------------------------------------------
\item not harmonic, i.e., $div R \ne 0$ and moreover $div C \ne 0$, $div P \ne 0$.
%-------------------------------------------------
\end{enumerate}
\end{rem}
%=============================================================================
\begin{rem}
In \cite{DD91a} Defever and Deszcz showed that if $Q(S,R) = 0$, then $R = L S\wedge S$ for some scalar $L$ if $S$ is not of rank 1. Recently, Shaikh and Kundu \cite{SKppsn} presented a metric (Example 1, \cite{SKppsn}) with $Q(S, R) = 0$, on which $S$ is not of rank 1 and $R = e^{x^1} S\wedge S$. It is interesting to mention that the rank of the Ricci tensor of the pp-wave metric is 1 and here $R\ne 0$ but $S\wedge S = 0$.
\end{rem}
\begin{rem}
It is well-known that every Ricci recurrent space with $\Pi$ as the 1-form of recurrency, is weakly Ricci symmetric with solution $(\Pi, 0, 0)$. It is interesting to mention that there may infinitely many solutions for a weakly Ricci symmetric manifold. The pp-wave metric (\eqref{gppwm} with condition \eqref{cond}) is weakly Ricci symmetric with solution
$$\Pi = \left\{\Pi_1, 0, \frac{f \left(h_{333}+h_{344}\right)-f_3 \left(h_{33}+h_{44}\right)}{f \left(h_{33}+h_{44}\right)}, \frac{f \left(h_{334}+h_{444}\right)-f_4 \left(h_{33}+h_{44}\right)}{f \left(h_{33}+h_{44}\right)}\right\},$$
$$\Omega = \left\{\Omega_1,0,0,0\right\} \ \mbox{ and}$$
$$\Theta = \left\{\frac{h_{133}+h_{144}}{h_{33}+h_{44}}-\Pi_1- \Omega_1, 0, 0, 0\right\},$$
where $\Pi_1$ and $\Omega_1$ are arbitrary scalars.
\end{rem}
%%%%%%%%%%%%%%%%%%%%%%%%%%%%%%%%%%%%%%%%%%%%%%%%%%%%%%%%%%%%%%%%%%%%%%%%%%%%%%%%%%%%%%%%%%%%%%%%%%%%%%%%%%%%%%%
%%%%%%%%%%%%%%%%%%%%%%%%%%%%%%%%%%%%%%%%%%%%%%%%%%%%%%%%%%%%%%%%%%%%%%%%%%%%%%%%%%%%%%%%%%%%%%%%%%%%%%%%
\indent Again it is clear that the pp-wave metric \eqref{pp-bc} in Brinkmann coordinates, is a special case of \eqref{gppwm} for $f\equiv -2$ and $h = -\frac{1}{2}H$ and hence satisfies \eqref{cond}. Therefore from Theorem \ref{mt}, we can state the following about the geometric properties of the metic \eqref{pp-bc}.
%=====================================================================================================
\begin{cor}\label{cor1}
The metric given in \eqref{pp-bc} possesses the following curvature properties:
\begin{enumerate}[label=(\roman*)]
%---------------------------------------------------------------------------------------------------
		\item $\kappa = 0$ and hence $R=W$ and $C=K$.
		%---------------------------------------------------------------------------------------------------
    \item $R$-space and $C$-space by Venzi for $\{1,0,0,0\}$. Hence from second Bianchi identity, the curvature 2-forms $\Omega_{(R)l}^m$ are recurrent for the 1-form $\{1,0,0,0\}$.
		%---------------------------------------------------------------------------------------------------
		\item Semisymmetric and hence Ricci semisymmetric, conformally semisymmetric and projectively semisymmetric. 
		%--------------------------------------------------------------------------------------------------
		\item If $4H_{34}^2+ \left(H_{33}-H_{44}\right)^2 \ne 0$, then its conformal curvature 2-forms $\Omega_{(C)l}^m$ are recurrent with 1-form of recurrency $\Pi$, given by
		$$\Pi_1 = 1, \ \ \Pi_2 = 0,$$
		$$\Pi_3 = \frac{2 H_{34} \left(H_{334}+H_{444}\right)+\left(H_{33}-H_{44}\right) \left(H_{333}+H_{344}\right)}{4H_{34}^2+ \left(H_{33}-H_{44}\right)^2},$$
		$$\Pi_4 = \frac{2 H_{34} \left(H_{333} + H_{344}\right) - \left(H_{33}-H_{44}\right) \left(H_{334}+ H_{444}\right)}{4H_{34}^2+ \left(H_{33}-H_{44}\right)^2}.$$
		%---------------------------------------------------------------------------------------------------
		\item If $H_{33}+H_{44} \ne 0$, then it is Ricci recurrent  with 1-form of recurrency $\Pi$, given by
		$$\Pi = \left\{1, 0, \frac{H_{333}+H_{344}}{H_{33}+H_{44}}, \frac{H_{334}+H_{444}}{H_{33}+H_{44}}\right\}.$$
		%----------------------------------------------------------------------------------------------------
		\item If $H_{33}+H_{44} \ne 0$, then it is weakly Ricci symmetric  with solution $(\Pi, \Omega, \Theta)$, given by
		$$\Pi = \left\{\Pi_1, 0, \frac{H_{333}+H_{344}}{H_{33}+H_{44}}, \frac{H_{334}+H_{444}}{H_{33}+H_{44}}\right\},$$
		$$\Omega = \left\{\Omega_1, 0, 0, 0\right\} \ \mbox{and}$$
		$$\Theta = \left\{\frac{H_{133}+H_{144}}{H_{33}+H_{44}}-\Pi_1-\Omega_1, 0, 0, 0\right\},$$
		where $\Pi_1$ and $\Omega_1$ are arbitrary scalar.
		%----------------------------------------------------------------------------------------------------
		\item If $H_{33}+H_{44} \ne 0$, then it is weakly cyclic Ricci symmetric  with solution $(\Pi, \Omega, \Theta)$, given by
		$$\Pi = \left\{\Pi_1, 0, \frac{H_{333}+H_{344}}{H_{33}+H_{44}}, \frac{H_{334}+H_{444}}{H_{33}+H_{44}}\right\},$$
		$$\Omega = \left\{\Omega_1, 0, \frac{H_{333}+H_{344}}{H_{33}+H_{44}}, \frac{H_{334}+H_{444}}{H_{33}+H_{44}}\right\} \ \mbox{and}$$
		$$\Theta = \left\{\frac{3 \left(H_{133}+H_{144}\right)}{H_{33}+H_{44}}-\Pi_1-\Omega_1, 0, \frac{H_{333}+H_{344}}{H_{33}+H_{44}}, \frac{H_{334}+H_{444}}{H_{33}+H_{44}}\right\},$$
		where $\Pi_1$ and $\Omega_1$ are arbitrary scalar.
		%----------------------------------------------------------------------------------------------------
		\item Ricci simple \cite{MS16} (i.e., $S = \alpha \eta\otimes\eta$) for
		$$\alpha = \frac{1}{2}\left(H_{33}+H_{44}\right) \ \ \mbox{ and } \ \ \eta = \left\{1, 0, 0, 0\right\}$$
		and hence $S\wedge S =0$ and $S^2 = 0$. Here $||\eta|| = 0$ and $\nabla \eta = 0$.
		%-----------------------------------------------------------------------------------------------------
		\item $Q(S, R) = Q(S, C) =0$ but $R$ or $C$ is not a scalar multiple of $S\wedge S$ as $S$ is of rank 1.
		%------------------------------------------------------------------------------------------------------
		\item $C\cdot R =0$ and hence $C\cdot S =0$, $C\cdot C =0$ and $C\cdot P =0$.
		%------------------------------------------------------------------------------------------------------
		\item $P\cdot R = 0$ but $P\cdot \mathcal R \ne 0$. Also but $P\cdot S =P\cdot \mathcal S =0$.
		%------------------------------------------------------------------------------------------------------
		\item Ricci tensor is Riemann compatible as well as Weyl compatible.
		%------------------------------------------------------------------------------------------------------
		\item $P\cdot P = -\frac{1}{3}Q(S,P)$.
		%------------------------------------------------------------------------------------------------------
\end{enumerate}
\end{cor}
%%%%%%%%%%%%%%%%%%%%%%%%%%%%%%%%%%%%%%%%%%%%%%%%%%%%%%%%%%%%%%%%%%%%%%%%%%%%%%%%%%%%%%%%%%%%%%%%%%%%%%%%%%%%%%%
\indent Again the non-vacuum pp-wave solution presented by Sippel and Goenner \cite{SG86} is a special case of \eqref{pp-bc} for $H(x,x^3,x^4) = 2 a_1 e^{a_2 x^3 - a_3 x^4}$. Hence the line element is explicitly given by:
\be\label{gppsg}
ds^2 = 2 a_1 e^{a_2 x^3 - a_3 x^4} (dx)^2 + 2 dx dr+[(dx^3)^2+(dx^4)^2].
\ee
Now the geometric properties of the metric \eqref{gppsg} can be stated as follows: 
%=====================================================================================================
\begin{cor}\label{cor2}
The metric given in \eqref{gppsg} possesses the following curvature properties:
\begin{enumerate}[label=(\roman*)]
%---------------------------------------------------------------------------------------------------
		\item it satisfies the curvature conditions (i)-(xiii) of Corollary \ref{cor1} with different associated 1-forms of the corresponding structures,
		%---------------------------------------------------------------------------------------------------
		\item moreover it is recurrent for the 1-form of recurrency $\{0,0,2 a_2,-2 a_3\}$. Hence it is Ricci recurrent, conformally recurrent and projectively recurrent. Also it is semisymmetric and hence Ricci semisymmetric, conformally semisymmetric and projectively semisymmetric.
\end{enumerate}
\end{cor}
%%%%%%%%%%%%%%%%%%%%%%%%%%%%%%%%%%%%%%%%%%%%%%%%%%%%%%%%%%%%%%%%%%%%%%%%%%%%%%%%%%%%%%%%%%%%%%%%%%%%%%%%%%%%%%%
\indent Again the generalized plane wave metric \cite{ste03} is given by
\be\label{gpwm}
ds^2 = 2 H(x,x^3,x^4) (dx)^2 + 2 dx dr+(dx^3)^2+(dx^4)^2,
\ee
where $H(x,x^3,x^4) = a_1 (x^3)^2 + a_2 (x^4)^2 + a_3 x^3 x^4 + a_4 x^3 + a_5 x^4 + a_6$, $a_i$'s are scalar. Hence it is a special case of \eqref{pp-bc} and we can state the following.
%=====================================================================================================
\begin{cor}\label{cor3}
The metric given in \eqref{gpwm} possesses the following curvature properties:
\begin{enumerate}[label=(\roman*)]
		%------------------------------------------------------------------------------------------------------
		\item it satisfies the curvature conditions (i)-(ix) of Corollary \ref{cor1} with different associated 1-forms of the corresponding structures,
		%------------------------------------------------------------------------------------------------------
		\item moreover it is locally symmetric and hence Ricci symmetric, conformally symmetric and projectively symmetric.
		%------------------------------------------------------------------------------------------------------
\end{enumerate}
\end{cor}
%============================================================================================
\indent From Corollary \ref{cor1},  Corollary \ref{cor2} and  Corollary \ref{cor3}, we can state the following about the recurrent structure on a semi-Riemannian manifold.
\begin{rem}
From Corollary \ref{cor2} we see that the metric \eqref{gppsg} is recurrent but not locally symmetric and from Corollary \ref{cor1} we see that the metric \eqref{pp-bc} is Ricci recurrent but not recurrent. These results support the well-known facts that every locally symmetric manifold is recurrent but not conversely, and every recurrent manifold is Ricci recurrent but not conversely.
\end{rem}
%
%%%%%%%%%%%%%%%%%%%%%%%%%%%%%%%%%%%%%%%%%%%%%%%%%%%%%%%%%%%%%%%%%%%%%%%%%%%%%%%%%%%%%%%%%%%%%%%%%%%%%%%%%%%%%%%%
%                                                  Energy-momentum tensor
%%%%%%%%%%%%%%%%%%%%%%%%%%%%%%%%%%%%%%%%%%%%%%%%%%%%%%%%%%%%%%%%%%%%%%%%%%%%%%%%%%%%%%%%%%%%%%%%%%%%%%%%%%%%%%%%
\section{\bf Energy-momentum tensor of generalized pp-wave metric}\label{E-M}
%%%%%%%%%%%%%%%%%%%%%%%%%%%%%%%%%%%%%%%%%%%%%%%%%%%%%%%%%%%%%%%%%%%%%%%%%%%%%
In this section we discuss about the energy-momentum tensor of the generalized pp-wave metric \eqref{gppwm} and also the other special forms, such as \eqref{pp-bc}, \eqref{gppsg} and \eqref{gpwm}. From the values of the energy momentum tensor $T$ of the generalized pp-wave metric \eqref{gppwm}, we can conclude that $T$ is of rank 1 if the cosmological constant is zero and \eqref{cond} holds. In this case 
$$T = \frac{c^4 \left(h_{33} + h_{44}\right)}{4 \pi  f G} \eta\otimes\eta, \ \ \eta=\{1,0,0,0\}.$$
Again it is easy to check that $\eta$ is null. Thus we can state the following:
\begin{thm}
For zero cosmological constant, the generalized pp-wave metric \eqref{gppwm} is a pure radiation metric if and only if it is a pp-wave metric.
\end{thm}
\begin{cor}
For zero cosmological constant, the pp-wave metric (\eqref{gppwm} with \eqref{cond} or \eqref{pp-bc}) is a pure radiation metric.
\end{cor}
%=======================================
Again from the values of the components of $\nabla T$ of the metric \eqref{gppwm}, we get
$$T_{11,1} + T_{11,1} + T_{11,1}=\frac{3 c^4 \left(h_{133}+h_{144}\right)}{4 \pi  f G},$$
\beb
T_{11,3} + T_{13,1} + T_{31,1}&=&\frac{c^4}{4 \pi  f^4 G} \left(f^3 h_{333}+f^3 h_{344}-f_{333} f^2 h\right.\\
&&-f_{344} f^2 h+f_{33} f^2 h_3+f_{44} f^2 h_3-f_3 f^2 h_{33}-f_3 f^2 h_{44}+4 f_3 f_{33} f h\\
&&\left.+2 f_4 f_{34} f h+2 f_3 f_{44} f h-f_3^2 f h_3-f_4^2 f h_3-3 f_3^3 h-3 f_3 f_4^2 h\right),
\eeb
\beb
T_{11,4} + T_{14,1} + T_{41,1}&=&\frac{c^4}{4 \pi  f^4 G} \left(f^3 h_{334}+f^3 h_{444}-f_{334} f^2 h-f_{444} f^2 h\right.\\
&&+f_{44} f^2 h_4-f_4 f^2 h_{33}-f_4 f^2 h_{44}+2 f_4 f_{33} f h+2 f_3 f_{34} f h+f_{33} f^2 h_4\\
&&\left.+4 f_4 f_{44} f h-f_3^2 f h_4-f_4^2 f h_4-3 f_4^3 h-3 f_3^2 f_4 h\right),
\eeb
$$T_{12,3} + T_{23,1} + T_{31,2}=\frac{c^4 \left(f^2 f_{333}+f^2 f_{344}+3 f_3^3+3 f_4^2 f_3-4 f f_{33} f_3-2 f f_{44} f_3-2 f f_4 f_{34}\right)}{8 \pi  f^4 G},$$
$$T_{12,4} + T_{24,1} + T_{41,2}=\frac{c^4 \left(f^2 f_{334}+f^2 f_{444}+3 f_4^3+3 f_3^2 f_4-2 f f_{33} f_4-4 f f_{44} f_4-2 f f_3 f_{34}\right)}{8 \pi  f^4 G}$$
and
\beb
T_{11,3} - T_{13,1} &=& \frac{c^4}{8 \pi  f^4 G} \left(2 f^3 h_{333}+2 f^3 h_{344}-2 f_{333} f^2 h-2 f_{344} f^2 h\right.\\
&&-f_{33} f^2 h_3-f_{44} f^2 h_3-2 f_3 f^2 h_{33}-2 f_3 f^2 h_{44}+8 f_3 f_{33} f h+4 f_4 f_{34} f h\\
&&\left.+4 f_3 f_{44} f h+f_3^2 f h_3+f_4^2 f h_3-6 f_3^3 h-6 f_3 f_4^2 h\right),
\eeb
\beb
T_{11,4} - T_{14,1} &=& \frac{c^4}{8 \pi  f^4 G} \left(2 f^3 h_{334}+2 f^3 h_{444}-2 f_{334} f^2 h-2 f_{444} f^2 h-f_{33} f^2 h_4\right.\\
&&-f_{44} f^2 h_4-2 f_4 f^2 h_{33}-2 f_4 f^2 h_{44}+4 f_4 f_{33} f h+4 f_3 f_{34} f h+8 f_4 f_{44} f h\\
&&\left.+f_3^2 f h_4+f_4^2 f h_4-6 f_4^3 h-6 f_3^2 f_4 h\right),
\eeb
$$T_{12,3} - T_{13,2} = T_{23,1} - T_{21,3} = -\frac{c^4 \left(f^2 f_{333}+f^2 f_{344}+3 f_3^3+3 f_4^2 f_3-4 f f_{33} f_3-2 f f_{44} f_3-2 f f_4 f_{34}\right)}{8 \pi  f^4 G},$$
$$T_{12,4} - T_{14,2} = T_{24,1} - T_{21,4} = -\frac{c^4 \left(f^2 f_{334}+f^2 f_{444}+3 f_4^3+3 f_3^2 f_4-2 f f_{33} f_4-4 f f_{44} f_4-2 f f_3 f_{34}\right)}{8 \pi  f^4 G}.$$
Now putting the condition \eqref{cond} to above, we get
\be\label{codtcom}
\left.\begin{array}{l}
T_{11,3} - T_{13,1} = \frac{c^4 \left(f h_{344}+f h_{333}-f_3 h_{33}-f_3 h_{44}\right)}{4 \pi  f^2 G},\\
$$T_{11,4} - T_{14,1} = \frac{c^4 \left(f h_{444}+f h_{334}-f_4 h_{33}-f_4 h_{44}\right)}{4 \pi  f^2 G},
\end{array}\right\}
\ee
and
\be\label{cstcom}
\left.\begin{array}{l}
T_{11,1} + T_{11,1} + T_{11,1}=\frac{3 c^4 \left(h_{144}+h_{133}\right)}{4 \pi  f G},\\
T_{11,3} + T_{13,1} + T_{31,1}=\frac{c^4 \left(f h_{344}+f h_{333}-f_3 h_{33}-f_3 h_{44}\right)}{4 \pi  f^2 G},\\
T_{14,1} + T_{41,1} + T_{11,4}=\frac{c^4 \left(f h_{444}+f h_{334}-f_4 h_{33}-f_4 h_{44}\right)}{4 \pi  f^2 G}.
\end{array}\right\}
\ee
%======================================================
Now from \eqref{dtcom} and \eqref{cstcom} we can state the following:
\begin{pr}
The energy momentum tensor $T$ of the pp-wave metric \eqref{gppwm} with the condition \eqref{cond} is\\
i) parallel if
$$h_{144}+h_{133} = f h_{344}+f h_{333}-f_3 h_{33}-f_3 h_{44} = f h_{444}+f h_{334}-f_4 h_{33}-f_4 h_{44} = 0,$$
(ii) Codazzi type if
$$f h_{344}+f h_{333}-f_3 h_{33}-f_3 h_{44} = f h_{444}+f h_{334}-f_4 h_{33}-f_4 h_{44} = 0,$$
(iii) cyclic parallel if $$h_{144}+h_{133} = f h_{344}+f h_{333}-f_3 h_{33}-f_3 h_{44} = f h_{444}+f h_{334}-f_4 h_{33}-f_4 h_{44} = 0.$$
\end{pr}
\begin{pr}
The energy-momentum tensor of the generalized pp-wave metric of the form \eqref{pp-bc} is
(i) parallel if
$$h_{144}+h_{133} = h_{344}+h_{333} = h_{444}+h_{334} = 0,$$
(ii) Codazzi type if
$$h_{344}+h_{333} = h_{444}+h_{334} = 0,$$
(iii) cyclic parallel if $$h_{144}+h_{133} = h_{344}+h_{333} = h_{444}+h_{334} = 0.$$
\end{pr}
%=================================================================================
Now we can conclude the following:
\begin{thm}
On a pp-wave spacetime (endowed with the metric \eqref{gppwm} with \eqref{cond} or with the metric \eqref{pp-bc})\\
(i) the energy-momentum tensor is parallel if and only if it is cyclic parallel,\\
(ii) the energy-momentum tensor is Codazzi type if it is cyclic parallel but not conversely (see Example \ref{ex1}),\\
(iii) the Ricci tensor is zero, i.e., the space is vacuum if $h$ is harmonic in $x^3$ and $x^4$, i.e., $h_{33} + h_{44} = 0$ and in this case $\nabla T = 0$.
\end{thm}
%=======================================================================================
\begin{exm}\label{ex1}
We now consider a special form of the generalized pp-wave metric \eqref{gppwm} as 
\be\label{spec}
ds^2=-2 e^{x+x^3-x^4} (dx)^2 + 2 dx dr - \frac{1}{2}e^{x^3-x^4}[(dx^3)^2+(dx^4)^2].
\ee
Then the non-zero components of its $R$, $\nabla R$, $S$ and $\nabla S$ are given by
$$R_{1313} = R_{1414} = e^{x+x^3-x^4}, \ \ \ R_{1313,1}= R_{1414,1}=e^{x+x^3-x^4},$$
$$S_{11}=4 e^x, \ \ \ S_{11,1}=4 e^x.$$
It is easy to check that the scalar curvature of this metric is zero and it is conformally flat. Now the non-zero components of its energy momentum tensor $T$ and its derivative $\nabla T$ are given by
$$T_{11}=\frac{c^4 e^{x-x^4} \left(2 e^{x^4}-\Lambda  e^{x^3}\right)}{4\pi G}, \ \ T_{12}=\frac{c^4 \Lambda }{8\pi G}, \ \ T_{33}=T_{44}=-\frac{c^4 \Lambda  e^{x^3-x^4}}{16\pi G},$$
$$T_{11,1}=\frac{c^4 e^x}{2\pi G}.$$
Thus we can easily check that the Ricci tensor and the energy momentum tensor of \eqref{spec} are codazzi type but not cyclic parallel.
\end{exm}
%
%%%%%%%%%%%%%%%%%%%%%%%%%%%%%%%%%%%%%%%%%%%%%%%%%%%%%%%%%%%%%%%%%%%%%%%%%%%%%%%%%%%%%%%%%%%%%%%%%%%%%%%%%%%%%%%%
%                                                  Comparison
%%%%%%%%%%%%%%%%%%%%%%%%%%%%%%%%%%%%%%%%%%%%%%%%%%%%%%%%%%%%%%%%%%%%%%%%%%%%%%%%%%%%%%%%%%%%%%%%%%%%%%%%%%%%%%%%
\section{\bf Robinson-Trautman metric and generalized pp-wave metric}\label{compar}
%%%%%%%%%%%%%%%%%%%%%%%%%%%%%%%%%%%%%%%%%%%%%%%%%%%%%%%%%%%%%%%%%%%%%%%%%%%%%%%%%%%
Recently, Shaikh et al. \cite{SAA17} studied the curvature properties of Robinson-Trautman metric. The line element of Robinson-Trautman metric in $\{t,r,x^3,x^4\}$-coordinate is given by
\be\label{rtm}
ds^2 = -2(a - 2 b r - q r^{-1})dt^2 + 2 dt dr - \frac{r^2}{f^2} [(dx^3)^2+(dx^4)^2],
\ee
where $a, b, q$ are constants and $f$ is a function of the real variables $x^3$ and $x^4$. In this section we make a comparison between the curvature properties of the Robinson-Trautman metric \eqref{rtm} and generalized pp-wave metric \eqref{gppwm} as well as the pp-wave metric \eqref{pp-bc}.
%=====================================================
\begin{thm}
The Robinson-Trautman metric \eqref{rtm} and generalized pp-wave metric \eqref{gppwm} have the following similarities and dissimilarities: \\
\textbf{A. Similarities:}\\
(i) Both the metrics are 2-quasi-Einstein,\\
(ii) both are generalized quasi-Einstein in the sense of Chaki,\\
(iii) Ricci tensors of both the metrics are Riemann compatible as well as Weyl compatible,\\
\textbf{B. Dissimilarities:}\\
(iv) \eqref{rtm} is Deszcz pseudosymmetric whereas \eqref{gppwm} is Ricci generalized pseudosymmetric,\\
(v) the conformal curvature 2-forms are recurrent for \eqref{rtm} but not recurrent for \eqref{gppwm},\\
(vi) the metric \eqref{rtm} is Roter type and hence $Ein(2)$ but \eqref{gppwm} is not Roter type but $Ein(3)$.
\end{thm}
%======================================================
\begin{thm}
The Robinson-Trautman metric \eqref{rtm} and pp-wave metric \eqref{pp-bc} have the following similarities and dissimilarities: \\
\textbf{A. Similarities:}\\
(i) for both the metrics, the conformal curvature 2-forms are recurrent,\\
(ii) Ricci tensor of both the metrics are Riemann compatible as well as Weyl compatible,\\
\textbf{B. Dissimilarities:}\\
(iii) the metric \eqref{rtm} is 2-quasi-Einstein, where as \eqref{pp-bc} is Ricci simple and hence quasi-Einstein,\\
(iv) \eqref{rtm} is Deszcz pseudosymmetric whereas \eqref{pp-bc} is semisymmetric,\\
(v) \eqref{rtm} is pseudosymmetric due to conformal curvature tensor whereas \eqref{pp-bc} is semisymmetric due to conformal curvature tensor,\\
(vi) the metric \eqref{rtm} is $S\wedge S \ne 0$ but \eqref{gppwm} $Ein(3)$,\\
(vii) the metric \eqref{rtm} is Roter type and hence $Ein(2)$ but \eqref{gppwm} is not Roter type but $Ein(3)$.
\end{thm}
%%%%%%%%%%%%%%%%%%%%%%%%%%%%%%%%%%%%%%%%%%%%%%%%%%%%%%%%%%%%%%%%%%%%%%%%%%%%%%%%%%%%%%%%%%%%%%%%%%%%%%%%%%%%%%%%%
%%%%%%%%%%%%%%%%%%%%%%%%%%%%%%%%%%%%%%%%%%%%%%%%%%%%%%%%%%%%%%%%%%%%%%%%%%%%%%%%%%%%%%%%%%%%%%%%%%%%%%%%%%%%%%%%%

%%%%%%%%%%%%%%%%%%%%%

%%%%%%%%%%%%%%%%%%%%%%%%%%%%%%%%%%%%%%%%%%%%%%%%%%%%%%%%%%%%%%%%%%%%%%%%%%%%%%%%%%%%%%%%%%%%%%%%%%%%
\end{document}